\documentstyle[11pt]{article}
\newcommand{\e}{\equiv}

\title{\Large\bf On Fermat's marginal note: a suggestion }
\author{ {\sc Nico F. Benschop} \\[1ex]
  {\it AmSpade Research  ~(n.benschop@chello.nl)}~~~18 April 1998}
\date{Presented at Netherlands Math Conference nr.33 ~(NMC-33)~~ Univ. Twente}

\begin{document}
\maketitle

\begin{abstract}
A suggestion is put forward regarding a partial proof of $FLT$ case1, which
is elegant and simple enough to have caused Fermat's enthusiastic remark
in the margin of his Bachet edition of Diophantus' $Arithmetica$.
It is based on an extension of Fermat's Small Theorem ($FST$) to mod $p^k$
for any $k>$0, and the cubic roots of 1 mod $p^k$ for primes $p=1$ mod 6.
For this solution in residues the exponent $p$ distributes over a sum, which
blocks extension to equality for integers, providing a partial proof of
$FLT$ case1 for all $p=1$ mod 6. This simple solution begs the question why
it was not found earlier. Some mathematical, historical and psychological
reasons are presented.

In a companion paper, on the triplet structure of Arithmetic mod $p^k$, this
cubic root solution is extended to the general rootform of $FLT$ mod $p^k$
(case1), called $triplet$. While the cubic root solution involves one inverse
pair: $a+a^{-1} \e -1$ mod $p^k ~(a^3 \e 1$ mod $p^k$), a triplet has three
inverse pairs in a 3-loop: $a+b^{-1} \e b+c^{-1} \e c+a^{-1} \e -1$ mod $p^k$
where $abc\e 1$ mod $p^k$, which is not restricted to $p$-th power residues
(for some $p \geq 59$) but applies to all residues in the group $G_k(.)$ of
units in the semigroup of multiplication mod $p^k$.
\end{abstract}

\section{ Introduction} 

Around 1637 Fermat discovered his Small Theorem ($FST$): $n^p \e n$ mod $p$ for
prime $p$ and all integers $n$, probably inspired by Pascal's triangle: the
multiplicative (factorial) structure of the coefficients in the expansion
of $(a+b)^p$. Only if $p$ is prime does it divide the binomial coefficient
of each of the $p-1$ $mixed$ terms, that is: except $a^p$ and $b^p$. Hence
$p$ divides $(a+b)^p-(a^p+b^p)$, in other words $(a+b)^p \e a^p+b^p$ mod $p$,
so exponent $p$ distributes over a sum (mod $p$).

One wonders, as possibly Fermat did, if this equivalence could hold mod $p^k$
for $k>$1 and some special $a,b$ - thus extending $FST$ to higher precisions
$k$. It will be shown that a straightforward extension of $FST$ to mod $p^k$
plays, for $p$=1 mod 6, an essential role in a special solution of normalized
form $FLT$ mod $p^k$: ~$a^p+b^p \e -1$ mod $p^k$. Here exponent $p$ distributes
over a sum, yielding a partial proof of the $FLT$ inequality for integers
(in case1: $a,b$ coprime to $p$).

\subsection{ Extending $FST$ to mod $p^k$ for $k>$1} 

Notice that $n^p \e n$ mod $p$ implies $n^{p-1} \e 1$ mod $p$ for $n \neq 0$ mod $p$,
and in fact all $p-1$ non-zero residues mod $p$ are known to form under
multiplication a cyclic group of order $p-1$. There are $p^{k-1}$ multiples of
$p$ among the $p^k$ residues mod $p^k$.
So $(p-1)p^{k-1}$ residues are coprime to $p$. They form the group $G_k$ of
$units$ in the semigroup $Z_k(.)$ of multiplication mod $p^k$. For each
$k>0$ there is necessarily a cyclic subgroup of order $p-1$, called the
{\bf core} $A_k$ of $G_k$. Clearly each $n$ in core has $n^p  \e n$ mod $p^k$.

Actually, units group $G_k$ is known to be cyclic for $p>$2 and all $k>$0.
Its order, as product of two coprime factors, implies it is a direct product
$G_k  \e A_k.B_k$ of two subgroups, namely {\bf core} ~$A_k$ of order $p-1$ and
{\bf extension} subgroup $B_k$ of order $p^{k-1}$. Each $n$ in core $A_k$
satisfies $n^p  \e n$ mod $p^k$, which clearly is a generalization of Fermat's
Small Theorem ($FST$) mod $p$ to mod $p^k$.

Furthermore, the choice of modulus $p^k$ yields every $p$-th iteration of
a generator of $G_k$, thus all $p$-th power residues, to form a subgroup
$F_k$ of order $|G_k|/p=(p-1).p^{k-2}$. For $k$=2 we have $A_2  \e F_2$: then
the core is the 'Fermat' subgroup of $p$-th power residues. In general we
have $|A_k|=p-1=|G_k|/p^{k-1}=|F_k|/p^{k-2}$, and $A_k=\{n^{p^{k-1}}\}$
mod $p^k$ for all $n$ in $G_k$.

\section{Solution of $FLT$ mod $p^k$ in Core has the $EDS$ property} 

The exponent distributes over a sum for an $FLT$ mod $p^k$ solution
{\bf in core}, because ~$(a+b)^p  \e a+b  \e a^p+b^p$
mod $p^k$, where $a^p  \e a$ and $b^p  \e b$ mod $p^k$. Such solution is said to
have the $EDS$ property: Exponent $p$ Distributes over a Sum.

If such a solution exists, as for each prime $p$=1 mod 6 (see further:
{\it cubic roots}), then  it cannot hold for integers, providing a
direct proof of integer inequality after all, despite Hensel's lemma of
infinite extension, described next. A solution $in~core$, having the
$EDS$ property, implies the $FLT$ (case1) inequality for integers.
Apart from a scaling factor, the cubic root solution is in fact [1]
the only one with all three terms in core $A_k ~(k \geq 3)$.

\subsection{ Hensel's extension lemma is no obstacle to a direct $FLT$ proof} 

Observe that for $k \geq 2$ core $A_k$ consists of $p$-th power residues.
The group of units is cyclic: $G_k  \e g^*$ with some generator $g$, and for
instance $|G_2|=(p-1)p$, so each $p$-th iteration of $g$ is in core $A_2$
which is a subgroup of order $p-1$.

It is easily verified that the two least significant digits of any $p$-ary
coded number determine if it is a $p$-th power residue, namely ~$iff$~
it is in core $A_2$. If so, then any more significant extension is also a
$p$-th power residue. This is known as {\it Hensel's extension lemma} (1913)
or the {\it Hensel lift}. This lemma implies that each $FLT$ mod $p^k$
solution is an more-significant digit ($msd$) extension of a solution mod $p^2$.

This lemma prevented the search for a direct $FLT$ proof via residues,
by the unwarranted conclusion that inequality for integers cannot be derived
from equivalence mod $p^k$. In fact, the solutions of $FLT$ mod $p^k$ (case1)
can all be shown to have exponent $p$ distributing over a sum, the "$EDS$"
property (or a variation of it) [1, lem3.1]~ yielding inequality for integer
$p$-th powers $<p^{pk}$.

\subsection{ Cubic roots of 1 mod $p^k$ ~sum to 0 mod $p^k$} 

Additive analysis shows that $each$ core subgroup $S \supset 1$, hence of
order $|S|$ dividing $p-1$, sums to 0 mod $p^k$ ($core~theorem$). If 3 divides
$p-1$, hence $p$=1 mod 6, the subgroup $S=\{a,~a^2,~a^3=1\}$ of the three
cubic roots of 1 mod $p^k$ sum to 0 mod $p^k$ , solving $FLT$ mod $p^k$.

For $|S|=3$ this zero sum is easily derived by simple means, without the
elementary semigroup concepts necessary to derive the additive $core~thm$ in
general. So Fermat might have derived this cubic root solution of $FLT$ mod
$p^k$ for $p$=1 mod 6, starting at $p$=7.~~ A simple {\bf proof} of the cubic
roots of 1 mod $p^k$ to have zero sum follows now, showing $a+b \e -1$ mod
$p^k$ to coincide with $ab  \e 1$ mod $p^k$.~~ Notice ~$a+b=-1$ ~to yield
~$a^2+b^2=(a+b)^2-2ab=1-2ab$, ~and:

~~~~ $a^3+b^3=(a+b)^3-3(a+b)ab=-1+3ab$. ~~The combined sum is $ab-1$:

~~~~ $\sum_{i=1}^3(a^i+b^i)=\sum_{i=1}^3 a^i + \sum_{i=1}^3 b^i= ab-1$
~mod $p^k$. ~Find ~$a,b$~ for ~$ab \e 1$ mod $p^k$.

Since $n^2+n+1=(n^3-1)/(n-1)$=0 ~for~ $n^3$=1 ($n \neq 1$), we have~ $ab$=1
mod $p^{k>0}$ ~if~ $a^3 \e b^3\e 1$ mod $p^k$, so 3 must divide $p-1~(p$=1 mod 6).

\subsection{ Proof of $FLT ~(case1)$ for $p=3, ~~5, ~~7$} 

Consider now only $FLT$ case1.
As mentioned earlier, the known Hensel extension lemma yields each solution
of $FLT$ mod $p^k$ to be an extension of a solution mod $p^2$, so analysis of
the normed $a^p+b^p \e -1$ mod $p^2$ is necessary and sufficient for the existence
of solutions at $p$ for any $k$.

For $p$=3 we have $|G_2|$=2.3 with core $A_2=\{-1,~1\}$, so core-pairsums
yield $\{-2,~0,~2\}$ which are not in core $A_2$, hence are not $p$-th power
residues. So the $FLT$ inequality holds for $p$=3.

For $p$=5: $G_2=2^*$ with ~$|G_2|$=4.5, and core $A_2=(2^5)^*=7^*$ mod 25. So
$A_2 \e \{7,-1,-7,1\}$ and non-zero coresums $\pm \{2,6,8,14\}$ which are not in
core $A_2$, hence are not $p$-th power residues, and thus $FLT$ holds for $p$=5.

For $p$=7: $G_2=3^*$ (order 6.7=42) and core $A_2=(3^7)^*=43^*=
\{43,42,66,24,25,01\}$ (base 7). The sum of {\bf cubic roots} of 1:
$\{42,24,01\}$ yields equivalence mod $7^2$, which necessarily yields
inequality for integers due to the $EDS$ property. So for $p$=7, and in fact
for all $p\e 1$ mod 6, $FLT$ (case1) holds for the corresponding cubic root
solutions.

\section{ Triplets as general root-form of $FLT$ mod $p^k$} 

A cubic root solution involves one inverse pair:
~~$a+a^{-1} \e -1$ mod $p^k ~(a^3\e 1$,~$a \neq 1, ~a^{-1} \e a^2)$.
The question remains if possibly other solutions to $FLT$ mod $p^k$ exist,
which can be answered by elementary semigroup techniques. In fact there is
precisely $one$ other solution type involving $three$ inverse pairs in a
successor coupled $loop$ of length 3, called

~~~{\bf triplet}:
    ~~$a+1 \e -b^{-1}, ~~b+1 \e -c^{-1}, ~~c+1 \e -a^{-1}$ mod $p^k$,
    ~where $abc\e 1$ mod $p^k$.

If $a \e b \e c$ then this reduces to the cubic root solution, which holds
for each prime $p\e 1$ mod 6. Triplet solutions occur for some primes
$p \geq 59$. A variant of the $EDS$ property can be derived for them
[1, lem3.1] sothat $FLT$ (case1) holds for all primes $p>$2.

\section{Summary} 

{\bf Did Fermat find the cubic root solution?} ~~If Fermat knew the cubic
root solution for $p\e 1$ mod 6, and also could prove it by elementary
means, as shown earlier, this might explain his enthusiastic note in the
margin, about a beautiful proof.

However, to complete the proof of case1 it is required to show that the cubic
root solutions are the $only$ solution type, which they in fact are not.
It seems very unlikely that he knew about the triplets, which start at $p$=59.
So he probably let the problem rest, realizing the cubic roots are only
a partial proof of $FLT$ case1. Another obstacle would be $case_2$ where
$p$ divides one of $x,y,z$, which requires a somewhat different approach [1].

{\bf Experimenting with $p=3, ~5, ~7$ (mod $p^2$)}: ~On Fermat's conjecture
($FLT$) of the sum of two $p$-th powers never to yield a $p$-th power (for
$p>$2), consider the next assumption about what he might have discovered, 
using means available at that time (1640). As shown, $p=3$ and $p=5$ yield
no solution mod $p^2$, hence no solution of $FLT$ mod $p^k$ for any $k$.

However for $p=7$, using $p$-ary code for residues mod $p^k$ (prime $p>2,
~k$ digits) and experimenting mod $7^2$ ~($p=7,~k=2)$, it is readily verified
that~ $x^p+y^p \e z^p$ mod $p^2$ ~does have a solution with the cubic
root of unity: ~$a^3 \e 1$ mod $p^2$. ~In fact: ~$a+1 \e -a^{-1}$ ~mod $7^2
~(a \neq 1$ mod $7^2$),\\ or equivalently "one-complement" normal form:

~~~~$a+a^{-1} \e -1$  ~mod $p^2$,
   ~with $a \e 24$ (in 7-code, decimal 18) and $a^{-1} \e 42$ (decimal 30).

As shown, this cubic-root solution holds for every prime $p\e 1$ mod 6, and
moreover (and this is the clue): $a^p \e a$ mod $p^k$, for every $k>$0.
Because cubic root "$a$" is in a p-1 order $core$ subgroup of $p$-th power
residues in the units group $G_k$ mod $p^k$ in ring $Z_k ~(k>$1). So:

~~~~~  $a^p+a^{-p} \e a+a^{-1} \e -1 \e (-1)^p \e (a+a^{-1})^p$,
       ~~prime~ $p\e 1$ mod 6.

For this solution $in~core$ the Exponent $p$ Distributes over a Sum 
("$EDS$" property), which blocks extension to integer equality, proving
$FLT ~(case_1)$ for all such cubic root solutions.

\subsection{By 'modern' elementary concepts} 

Using elementary group concepts: the units group $G_k$ mod $p^k$ has order
$(p-1)p^{k-1}$, and $G_k$ is known to be cyclic for all $k>$0. The two
coprime factors imply $G_k$ to be a direct product of two cyclic groups

~~~~~~~~~ $G_k  \e  A_k.B_k$  with $core$ subgroup $|A_k|= p-1$,
          and $extension$ subgroup $|B_k| = p^{k-1}$.

Of course, these group theoretical arguments were not known in those days,
but the insight that $FLT$ mod $p^k$ does have a cubic root solution for every
prime $p\e 1$ mod 6 at each $k>$0, could very well be discovered by Fermat
-- first found by hand calculations mod $7^2$ (in 7-ary code for instance),
and then derived algebraically in general for all $p\e 1$ mod 6, as shown.

Actually, with the known group concepts as described above, the Bachet
margin might be large enough to sketch the essence of this proof for
the cubic root solutions, and the impossibility of extension to integer
equality (by the $EDS$ argument).

{\bf Cubic roots not the only solution form}.~
Assuming this solution occurred to Fermat, he must have realized that the
full proof (at least of case1: with $x,y,z$ coprime to p) would require to
show that the cubic root type of solution is the only solution type.

However this is not the case, as mentioned earlier: the general type of
solution is a "triplet": $a+b^{-1} \e b+c^{-1} \e c+a^{-1} \e -1$ mod $p^k$,
with $abc \e 1$ mod $p^k$, involving not one inverse-pair, but three
inverse-pairs in a loop of length 3 ~[1] - which only occurs for some
primes $p \geq 59$, and of which the cubic root solution is a special case.

A variation of the $EDS$ argument holds here, so again the $FLT$ inequality
for integers follows, proving $FLT$ case1 after showing that no other
solutions exist ( here the non-commutative function composition of semigroups
is essential, as applied to the two symmetries $-n$ and $n^{-1}$ of residue
arithmetic, as well as quadratic analysis mod $p^3$, see thm3.2 in [1] ).
\\Moreover, there is $FLT~case_2$~, with an approach as given in [1].

\section*{ Conclusion}

The above sketched cubic root solution is a partial proof of $FLT$ case1,
possibly known to Fermat. However, missing the triplets for $p \geq 59$
~(it is very doubtful that they could be found without computer experiments,
which are easy for present day PC's), and possibly lacking an approach for
$FLT$ case2, he probably was inclined to keep quiet about this partial
proof of $FLT$ case1.

It is surprising that the cubic root solution of $FLT$ mod $p^k$, and the
corresponding $EDS$ property of the exponent $p$ distributing over a sum,
did not surface long ago. Giving some thought to possible causes of this
delay, one might consider the following phenomena of mathematical,
historical and psychological nature.

\begin{enumerate}
\item
Missing link between $FST$ and $FLT$: ~It appears that, for unexplained
reasons, no link has been made between Fermat's Small and Last Theorem,
although both feature $p$-th powers: residues mod $p$ in the first, and
integers in the second. Furthermore, the $p-1$ cycle corresponding to
$n^p \e n$ is clearly common to the units group mod $p$ and mod $p^k$. So
it seems that the group structure, available since the second half of the
previous century, is not considered for some reason. Possibly because of
other promising approaches taken in the analysis of arithmetic
(e.g. Hensel's $p$-adic number theory, 1913).

\item
Dislike of exponentiation (~\^~) which is not associative, nor commutative,
nor does it disribute over addition (+). Closure properties holding for
(+,~.) do not hold for (~\^~). However, this situation is improved by taking
$p^k$ as modulus, because then the $p$-th power residues do form a subgroup
$F_k$ of the units group $G_k(.)$, which for $k$=2 in fact is the core of $G_2$
with a nice additive property (zero sum subgroups). None of these properties
is difficult to derive, and require only elementary semigroup concepts.
It seems that application of semigroups to arithmetic ran out of fashion, rather
being employed for the development of higher and more abstract purposes, such
as category theory. Clearly in the elliptic curve approach $modular~forms$,
which have good closure properties, are preferred over exponentiation
(Eichler [5]: "There are five basic arithmetic operations: addition,
subtraction, multiplication, division and modular forms").

\item The $notation$ under which $FST$ usually is known: $n^{p-1} \e 1$ mod $p$,
instead of the $fixed~point$ notation $n^p \e n$ mod $p$, is less like the $EDS$
form $(a+b)^p \e a^p+b^p$ mod $p$ which might have prompted Fermat to explore the
same mod $p^2$ and mod $p^k$ (2- and $k$- digit arithmetic), yielding the cubic
root solution mod $7^2$ as first example.

\item The Hensel lift. From Hensel's $p$-adic number theory (1913) derives the
known lemma about each solution of $FLT$ mod $p^k$ to be an extension of a
solution mod $p^2$. So analysis mod $p^2$ is necessary and sufficient for
root existence, with an $FLT_k$ root being a solution of $a^p+b^p \e -1$ mod $p^k$.
This infinite precision extension to all $k>$0, called the Hensel lift, is
the most often cited objection against a direct proof of $FLT$; an unwarranted
conclusion, due to the $EDS$ property of all $FLT_2$ roots: $F_2 \e A_2$ is core.
Moreover, analysis mod $p^3$ is necessary to describe the symmetries of $FLT_k$
roots, and characterize the general $triplet$ rootform [1]. \\
-- The irony is that the basis of the Hensel lift also carries the solution:
the $triplets$ follow from a detailed (computer-) analysis of the solutions of
$a^p+b^p \e -1$ mod $p^2$. Special attention was given to the role of inverse pairs,
as indicated by the cubic root solution, by using logarithmic code over a
primitive root of 1 mod $p^2$ ~(~[1] table 2~).

\item
Computer use. Simple computer experiments were very helpful, if not indispensible,
in discovering the triplets as general rootform -- given the importance of
inverse-pairs as evident from the cubic root solution. It is highly improbable
that they could be discovered any other way ($p \geq 59$). Unlike Pascal and
Leibniz, who even constructed their own accumulator and multiplier respectively,
the use of a computer seems to be avoided by some [5].

\item
The $FLT$ proof via elliptic curves (A.Wiles, {\it Annals of Mathematics},
May 1995) blocked further interest in a simpler proof, or at least rendered
such efforts irrelevant in the eyes of experts. The extreme complexity of that
proof seems to be interpreted as an advantage, rather than a disadvantage.

\item
The structure of a finite semigroup, starting with Schushkewitch's analysis
[4, appx.A] of its minimal ideal (1928), combined with the concept of
$rank$ to model the divisors of zero, is useful to derive the additive $core$
theorem in general ( [1]~thm1.1: each core subgroup $S \supset 1$ sums to
0 mod $p^k$ ). The clue of $FLT$ mod $p^k$ is to look for an additive property
in a multiplicative semigroup, although simple arithmetic suffices in the
cubic root case.

\item The role of Authority. An often used argument against a
simple $FLT$ proof is: So many eminent mathematicians have tried for so long
that it would have been found long ago. This argument of course does not take
into account the essential ingredients of the cubic-root \& triplet structure
of arithmetic mod $p^k$ [1] and the $EDS$ property, such as: semigroup
principles and computer experiments. These were only available since 1928 and
1950 respectively.\\ On the contrary, the lack of results in various directions
stresses the point that something was missing, requiring a different approach
(commutative arithmetic [+,~.~] in the context of associative function
composition: semigroups).

\item
Between disciplines: The application of semigroups, that is associative function
composition, to Arithmetic [4, p130] [6] . . . e.g. viewing its two symmetries:
complement $C(n)=-n$, inverse $I(n)=n^{-1}$, and the successor $S(n)=n+1$ as
functions [1, thm3.2], turned out to be a rare combination. The extreme
abstraction and specialisation in mathematics reduces the chances for a
serious consideration of an inter- disciplinary approach.
\end{enumerate}

{\Large\bf References} :

1. N. F. Benschop: " Triplets and Symmetries of Arithmetic mod $p^k$ "
\\ \hspace*{.5cm} http://www.iae.nl/users/benschop/nfb0.dvi
    ~and~ http://arXiv.org/abs/math.GM/0103014

--- The essence was presented at two international conferences: 2.(math), 3.(engineering):

2. ---: "On the semigroup of multiplication mod $p^k$, an extension of Fermat's
    Small Theorem, and \\. ~~~~its additive structure",
     {\it ~Semigroups and their Applications} ~(digest p7), ~Prague, July 1996.

3. ---: "On the additive structure of Multiplication mod $p^k$, Fermat's Small
      Theorem \\ \hspace*{.8cm} extended to $FLT$, and a new binary number code",
   \\ \hspace*{.5cm}
    {\it Logic and Architecture Synthesis} ~(digest pp.133-140) Grenoble, Dec 1996.

4. A. Clifford, G. Preston: {\it The Algebraic Theory of Semigroups},
   Vol I, pp.130-135,\\ \hspace*{.5cm} and Appendix A (p.207), ~AMS survey \#7, 1961.

5. A. Wiles ($BBC$ interview): ~http://www.bbc.co.uk/horizon/95-96/960115.html

6. S. Schwarz: "The Role of Semigroups in the Elementary Theory of Numbers", \\
   \hspace*{.5cm} Math.Slovaca V31 (1981) N4, pp.369-395.

\end{document}